\documentstyle[12pt, twoside]{article}
\input amssym.def
\input amssym.tex

\newcommand{\qed}{
  \strut\hfill
  \mbox{$\Box$} }

\newtheorem{theorem}{Theorem}[section]
\newtheorem{corollary}{Corollary}[section]
\newtheorem{lemma}{Lemma}[section]
\newtheorem{conjecture}{Conjecture}[section]

\newtheorem{example}{Example}[section]
\newtheorem{remark}{Remark}[section]
\newtheorem{definition}{Definition}[section]
\newtheorem{proposition}{Proposition}[section]
\newcommand{\1}{  \bf 1 }

\newcommand{\h}{ \frak h}
\newcommand{\D}{ \cal D }
\newcommand{\hD}{ \widehat {\cal D} }
\newcommand{\N}{ \Bbb N }
\newcommand{\Z}{ \Bbb Z }
\newcommand{\C}{ \Bbb C }

\def \be{\begin{equation}}
\def \ee{\end{equation}}
\def \bex{\begin{example}\label}
\def \bea{\begin{eqnarray} }
\def \eea{\end{eqnarray} }
\def \eex{\end{example}}

\newcommand{\W}{ {\cal W}_{1+\infty} }
\newcommand{\Wi}{ W_N}

\def \<{\langle}
\def \>{\rangle}

\def \a{\alpha }
\def \e{\epsilon }
\def \l{\lambda }

\def \b{\beta }

\def \gi {\gamma_i}
\def \be{\begin{equation}}
\def \ee{\end{equation}}
\def \bex{\begin{example}\label}
\def \bea{\begin{eqnarray} }
\def \eea{\end{eqnarray} }
\def \eex{\end{example}}
\def \bl{\begin{lem}\label}
\def \el{\end{lem}}
\def \bt{\begin{thm}\label}
\def \et{\end{thm}}
\def \bp{\begin{prop}\label}
\def \ep{\end{prop}}
\def \br{\begin{rem}\label}
\def \er{\end{rem}}
\def \bc{\begin{coro}\label}
\def \ec{\end{coro}}
\def \bd{\begin{de}\label}
\def \ed{\end{de}}

\newcommand{\la}{\langle}
\newcommand{\ra}{\rangle}

\def \l {\lambda}
\makeatletter \@addtoreset{equation}{section}

\makeatother \makeatletter
\newcommand{\bA}{{\bar A} }

\begin{document}
\title{ {
Representations of the vertex algebra $\W$ with a negative integer
central charge} }

\author{
  {Dra\v zen Adamovi\' c}
}
\date{}
\maketitle
\begin{abstract}
Let $\D$ be the Lie algebra of regular differential operators on
${\C} \setminus \{0\}$, and ${\hD}= {\D} + {\C} C$ be the central
extension of ${\D}$.
 Let $W_{1+\infty,-N}$ be the vertex algebra associated to the
 irreducible vacuum $\hD$--module with the central charge $c=-N$.
 We show that $W_{1+\infty,-N}$ is a subalgebra of the Heisenberg
 vertex algebra $M(1)$ with $2 N$ generators, and construct $2
 N$--dimensional family of irreducible $W_{1+\infty,-N}$--modules.
 Considering these modules as $\hD$--modules, we identify the
 corresponding highest weights.
\end{abstract}

\section{Introduction}

Let $\D$ be the Lie algebra of regular differential operators on
${\C} \setminus \{0\}$, and ${\hD}= {\D} + {\C} C$ be the central
extension of ${\D}$. In the representation theory of the Lie
algebra $\hD$ the most important are  the irreducible quasi-finite
highest weight modules. These modules were classified by V. Kac
and A. Radul in \cite{KR1}. In \cite{FKRW} was shown that the
language of vertex algebras is very useful in $\hD$--module
theory. In particular, on the irreducible vacuum $\hD$--module
$L(0,c,\hD)$ exists a natural structure of a simple vertex algebra
(cf. \cite{FKRW}, \cite{K}). This vertex algebra is usually
denoted by $W_{1+\infty,c}$. The results from \cite{KR1} give that
the representation theory of the vertex algebra $W_{1+\infty,c}$
is nontrivial only in the case $c \in {\Z}$. When $N \in {\N}$,
the irreducible modules for $W_{1+\infty,N}$ were classified in
\cite{FKRW}.

The structure of the vertex algebra $W_{1+\infty,-N}$ is much
complicated. The representation theory of $W_{1+\infty,-N}$ was
began by Kac and Radul in \cite{KR2}. They realized
$W_{1+\infty,-N}$ as a vertex subalgebra of the vertex algebra
$V_N$ constructed from Weyl algebra $W_N$ (we recall this result
in the Section \ref{reckr}). They also constructed a large class
of irreducible $W_{1+\infty,-N}$--modules.  The next step in this
direction was made by Wang in \cite{W1}, \cite{W2}. He considered
the case $c=-1$, and proved that $W_{1+\infty,-1}$ is isomorphic
to tensor product $W_{3,-2} \otimes H_0$, where $W_{3,-2}$ is a
vertex algebra associated with $W_3$--algebra with the central
charge $-2$ and $H_0$ is a Heisenberg vertex algebra. He also
classified all irreducible modules for $W_{3,-2}$ and
$W_{1+\infty,-1}$. The representations obtained in \cite{W2}
weren't identified as a highest weight $\hD$--modules.

In the present paper we will construct $2 N$--dimensional family
of irreducible $W_{1+\infty,-N}$--modules. Our family includes the
modules constructed in \cite{KR2}, and also coincides with the
modules from \cite{W2} in the case $c=-1$.

Let us explain our result in more details. We consider lattice
vertex superalgebra $V_L$ (cf. \ref{stwist}), and for suitably
chosen lattice $L$ we show that $V_N$ and $W_{1+\infty,-N}$ are
vertex subalgebras of $V_L$ (cf. Section \ref{wsec}). This fact is
in physical literature known as a bosonization of $\beta \gamma$
system (see \cite{FMS}, \cite{W1}). The embedding of
$W_{1+\infty,-N}$ into vertex superalgebra $V_L$ will imply that
$W_{1+\infty,-N}$ can be realized as a vertex subalgebra  of the
Heisenberg vertex algebra $M(1)$ with $2 N$--generators. We
explicitly identify the generators of $W_{1+\infty,-N}$ in terms
of Schur polynomials. Considering $M(1)$--modules $M(1,\l)$ we
obtain $W_{1+\infty,-N}$--modules $V(\l,-N)$ as a irreducible
subquotients of $M(1,\l)$ (see Section \ref{wsec}).  Considering
$V(\l,-N)$  as a  $\hD$--module,  we identify its highest weights.
As a result we get that all modules $V(\l,-N)$ are quasi-finite
$\hD$--modules.

\section{Vertex superalgebras}

In this section we recall the definition of vertex (super)algebra
and a few basic formulas. Fore more details about vertex
(super)algebras its representations,  and representation theory of
certain examples  of vertex (super)algebras see  \cite{B},
\cite{DL}, \cite{FLM}, \cite{K}, \cite{Li}, \cite{A}.

 Let $V$ be a vector space.  A
{\it field} is a series of the form

$$
  a (z) = \sum_{n \in {\Z } } a_n z^{- n - 1},
$$ where $a_n \in \mbox{End} V$ are such that for each $v \in V$
one has: $a_n v = 0 \quad \mbox{for} \quad n \gg 0$.  Here $z$ is
a formal indeterminate.

For a ${\Z} _2$-graded vector space $W=W^{even}+W^{odd}$ we write
$|u|\in {\Z}_2$, a degree of $u$, only for homogeneous elements:
$|u|=0$ for an even element $u\in W^0$ and $|u|=1$ for an odd
element $u\in W^1$. For any two $\Z _2$-homogeneous elements $u$
and $v$ we define $\e _{u,v}=(-1)^{|u||v|}\in {\Z}$.
\begin{definition}
 A {\it vertex superalgebra} is a quadruple
$(V,{\1},D,Y)$, where $V=V^{even}\oplus V^{odd}$ is a
$\Bbb{Z}_{2}$-graded vector space, $D$ is a
$\Bbb{Z}_{2}$-endomorphism of $V$, ${\1}$ is a specified vector
called the {\it vacuum} of $V$, and $Y$ is a linear map
\bea Y(\cdot,z):& &V\rightarrow ({\rm
End}V)[[z,z^{-1}]];\nonumber\\ & &a\mapsto Y(a,z)=\sum_{n\in
\Bbb{Z}}a_{n}z^{-n-1}\;\;(\mbox{where } a_{n}\in {\rm End}V)
\nonumber
 \eea
 such that
\bea (V1)& &\mbox{For any }a,b\in V, a_{n}b=0\;\;\;\mbox{ for }n
\mbox{ sufficiently large;}\nonumber \\ (V2)&
&[D,Y(a,z)]=Y(D(a),z)=\partial_z Y(a,z)\;\;\mbox{  for any }a\in
V;\nonumber
\\ (V3)& &Y({\1},z)=Id_{V}\;\;\;\mbox{(the identity operator of
$V$)};\nonumber \\ (V4)& &Y(a,z){\1}\in ({\rm End}V)[[z]] \mbox{
and }\lim_{z \rightarrow 0}Y(a,z){\1}=a\;\;\mbox{  for any }a\in
V;\nonumber \\ (V5)& &\mbox{For }\Bbb{Z}_{2}\mbox{ -homogeneous
}a,b\in V, \mbox{ the following {\it Jacobi identity} holds:}
\nonumber \eea
\bea &
&\;\;\;z_{0}^{-1}\delta\left(\frac{z_{1}-z_{2}}{z_{0}}\right)Y(a,z_{1})Y(b,z_{2})
-\varepsilon_{a,b}z_{0}^{-1}\delta\left(\frac{z_{2}-z_{1}}{-z_{0}}\right)
Y(b,z_{2})Y(a,z_{1})\nonumber \\ &
&=z_{2}^{-1}\delta\left(\frac{z_{1}-z_{0}}{z_{2}}\right)Y(Y(a,z_{0})b,z_{2}).\nonumber
\eea
\end{definition}
In the case $V= V^{even}$, i.e. when all vectors are even we say
that $V$ is vertex algebra.

 From the Jacobi identity follows
\bea \label{comut}
  \left[
    a_{m}, b_{n}
  \right] = \sum^\infty_{j = 0}
    { m \choose
      j} (a_{j} b)_{m + n - j}, \quad m, n \in {\Z},
\eea
\bea \label{5.2}
    Y (a_{- 1} b, z)  =
  : Y(a, z) Y(b, z):,
\eea
where the normally ordered product is defined, as usual, by
$$
  : Y (a, z) Y (b, z): = Y (a, z)_- Y(b, z) + Y (b, z)
  Y (a, z)_+
$$
and
$$
  Y(a, z)_+ = \sum_{j \in {\Z}_+} a_{j} z^{-j - 1},
  \quad Y(a, z)_- = Y(a, z) - Y (a, z)_+
$$
are the {\it annihilation\/} and the {\it creation\/} parts of $Y
(a, z)$.

\section{ The vertex  algebra $W_{1 + \infty, - N}$ }
\label{reckr}

In this section we recall some of the results from \cite{KR1},
\cite{KR2}, \cite{FKRW}.

Let $\D$ be the Lie algebra of complex  regular differential
operators on ${\C}^\times$ with the usual bracket, in an
indeterminate $t$.  The elements
$$
  J^l (k) = - t^{l + k} (\partial_t)^l \quad
  (k \in {\Z}, l \in {\Z}_+)
$$
form a basis of ${\D}$.  The Lie algebra ${\D}$ has the following
2-cocycle with values in ${\C}$ :

$$
  \Psi (f (t) \partial_t^m, g (t) \partial_t^n) =
  \frac{m! n!}{(m + n + 1)!}
  \mbox{Res}_{t = 0} f^{(n + 1)} (t) g^{(m)} (t) dt,
$$
where $f^{(m)} (t) = \partial^m_t f (t)$.  We denote by ${\widehat
{\D} } = {\D} \oplus {\C} C$, where $C$ is the central element,
the corresponding central extension of the Lie algebra ${\D}$.

Another important basis of ${\D}$ is
$$
  L^l (k) = - t^k D^l \quad
  (k \in {\Z}, l \in {\Z}_+)
$$
where $D = t \partial_t$.  These two bases are related by the
formula \cite{KR1}:
\bea \label{veza}
  J^l (k) = - t^k D (D - 1) \cdots (D - l + 1).
\eea

Given a sequence of complex numbers $\lambda = (\lambda_j)_{j
  \in {\Z}_+}$ and a complex number $c$ there exists a unique
irreducible ${\widehat {\D} }$-module $L (\lambda, c ; {\widehat
\D})$ which admits a non-zero vector $v_\lambda\/$ such that:
$$
  L^j (k) v_\lambda = 0 \quad \mbox{for} \quad
  k > 0,\; L^j (0) v_\lambda = \lambda_j v_\lambda, \; C = c I.
$$
This is called a highest weight module over $\widehat {\D}$ with
highest weight $\lambda\/$ and central charge $c$.  The module $L
(\lambda, c; \widehat \D)$ is called {\it quasifinite\/} if all
eigenspaces of ${D}$ are finite-dimensional (note that $D\/$ is
diagonalizable).  It was proved in (\cite{KR1} Theorem 4.2) that
$L (\lambda, c; \hat D)$ is a quasi-finite module if and only if
the generating series

$$ \Delta_\lambda (x) = \sum^\infty_{n = 0} \frac{x^n}{n!}
\lambda_n $$

has the form $$
  \Delta_\lambda (x) = \frac{\phi (x)}{e^x - 1},
$$ where

$$
  \phi (x) + c = \sum_i p_i (x) e^{r_i x} \quad
  \mbox{(a finite sum)},
$$

$p_i (x)$ are non-zero polynomials in $x\/$ such that $\sum_i p_i
(0) = c$ and $r_i\/$ are distinct complex numbers.  The numbers
$r_i\/$ are called the {\it exponents\/} of this module and the
polynomials $p_i (x)$ are called their {\it multiplicities}.

Recall now that the ${\hD} $-module $L(0, c; {\hat {\D} })$ has a
canonical structure of a vertex algebra with the vacuum vector
${\1} = v_0$ and generated by the fields $J^k (z) = \sum_{m \in
{\Z} } J^k (m) z^{-m - k - 1}$ \cite{FKRW}.

 In \cite{KR2}, the vertex algebra $L(0,-N, \D)$ was realized by
 using Weyl algebra ${\Wi}$ and the corresponding vertex algebra
 $V_N$.
The Weyl algebra ${\Wi}$ is an associative algebra over ${\C}$
generated by $a_i (m), a^{*}_i (m)$  $ (i=1,\dots, N; m \in {\Z})$
and $C$ with the following defining relations

$$ [a_i (m), a_j (n)]= [a^{*}_i (m) , a^{*}_j (n) ] = 0, \quad
[a^{*}_i (m), a_j (n) ] = \delta_{i,j} \delta_{m,-n} C
  $$
  for all $i,j \in \{ 1,\dots, N\}$, $n,m \in {\Z}$,
  and $C$ is central element.

  The vacuum ${\Wi}$--module $V_N$  is generated by one vector
  ${\1}$, and the following relations

$$ a_i (n) {\1} = 0, \quad n\ge 0; \quad
    a^{*}_i (n) {\1}= 0, \quad n> 0, \quad C {\1}= {\1}.$$

    Define the following fields acting on $V_N$.
 \bea \label{polja}
 a_i (z ) = \sum_{n \in {\Z} } a(n) z^{-n-1}, \quad
 a^{*} _i (z) = \sum_{n \in {\Z} }  a^{*} (n) z^{-n}.
 \eea

 Then there is a unique extensions of the fields (\ref{polja})
 such that $V_N$ becomes a vertex algebra (see \cite{KR2}).

 Denote by $W_{1 + \infty, -N}$ a vertex subalgebra of $V_N$
 generated by the fields
$$
  J^k (z) = - \sum^N _{i = 1} : a_i (z) \partial{^k} _{z} a^{*}
    _i (z): .
$$

\begin{proposition}(\cite{KR2}, Proposition 6.3)

We have an isomorphism of vertex algebras:
$$ L (0, -N; {\hat {\D}}) \simeq W_{1 + \infty, -N} $$
under which the fields (denoted by the same symbol) $J^k (z)$
correspond to each other.
\end{proposition}

\section{Lattice and Heisenberg vertex algebras } \label{stwist}

 Let $L$  be a  lattice.
  Set ${\h}={\C}\otimes_{\Z}L$ and
extend the ${\Z}$-form $ \la \cdot, \cdot \ra $ on $L$ to ${\h}$.
 Let $\hat{{\h}}={\C}[t,t^{-1}]\otimes {\h} \oplus {\C}c$ be the affinization of
${\h}.$ We also use the notation $h(n)=t^{n}\otimes h$ for $h\in
{\h}, n\in {\Z}$.

Set
$
\hat{{\h}}^{+}=t{\C}[t]\otimes
{\h};\;\;\hat{{\h}}^{-}=t^{-1}{\C}[t^{-1}]\otimes {\h}.
$
Then $\hat{{\h}}^{+}$ and $\hat{{\h}}^{-}$ are abelian subalgebras
of $\hat{{\h}}$. Let $U(\hat{{\h}}^{-})=S(\hat{{\h}}^{-})$ be the
universal enveloping algebra of $\hat{{\h}}^{-}$. Let ${\l} \in
{\h}$. Consider the induced $\hat{{\h}}$-module
\begin{eqnarray*}
M(1,{\l})=U(\hat{{\h}})\otimes _{U({\C}[t]\otimes {\h}\oplus
{\C}c)}{\C}_{\l}\simeq
S(\hat{{\h}}^{-})\;\;\mbox{(linearly)},\end{eqnarray*} where
$t{\C}[t]\otimes {\h}$ acts trivially on ${\C}$,
${\h}$ acting as $\la h, {\l} \ra$ for $h \in {\h}$
and $c$ acts on ${\C}$ as multiplication by 1. We shall write
$M(1)$ for $M(1,0)$.
 For $h\in {\h}$ and $n \in {\Z}$ write $h(n) =  t^{n} \otimes h$. Set
$
h(z)=\sum _{n\in {\Z}}h(n)z^{-n-1}.
$

Then $M(1)$ is vertex algebra which is generated by the fields
$h(z)$, $h \in {\h}$, and $M(1,{\l})$, for $\l \in {\h}$, are
irreducible modules for $M(1)$.

Let $\hat{L}$ be the canonical central extension of $L$ by the
cyclic group $\< \pm 1\>$:
\begin{eqnarray}\label{2.7}
1\;\rightarrow \< \pm 1\>\;\rightarrow \hat{L}\;\bar{\rightarrow}
L\;\rightarrow 1
\end{eqnarray}
with the commutator map $c(\alpha,\beta)=(-1)^{\< \alpha,\beta\>}$
for $\alpha,\beta \in L$. Let $e:  L\to \hat L$ be a section such
that $e_0=1$ and $\epsilon: L\times L\to \<\pm 1\>$ be the
corresponding 2-cocycle. Then
$\epsilon(\a,\b)\epsilon(\b,\a)=(-1)^{\<\a,\b\>},$
\begin{equation}\label{2c}
\e(\a,\b)\e(\a+\b,\gamma)=\e(\b,\gamma)\e(\a,\b+\gamma)
\end{equation}
and $e_{\a}e_{\b}= \e(\a,\b)e_{\a+\b}$ for $\a,\b,\gamma\in L.$
   Form the induced $\hat{L}$-module
\begin{eqnarray*}
{\C}\{L\}={\C}[\hat{L}]\otimes _{\< \pm 1\>}{\C}\simeq
{\C}[L]\;\;\mbox{(linearly)},\end{eqnarray*} where ${\C}[\cdot]$
denotes the group algebra and $-1$ acts on ${\C}$ as
multiplication by $-1$. For $a\in \hat{L}$, write $\iota (a)$ for
$a\otimes 1$ in ${\C}\{L\}$. Then the action of $\hat{L}$ on ${\C}
\{L\}$ is given by: $a\cdot \iota (b)=\iota (ab)$ and $(-1)\cdot
\iota (b)=-\iota (b)$ for $a,b\in \hat{L}$.

Furthermore we define an action of ${\h}$ on ${\C}\{L\}$ by:
$h\cdot \iota (a)=\< h,\bar{a}\> \iota (a)$ for $h\in {\h},a\in
\hat{L}$. Define $z^{h}\cdot \iota (a)=z^{\< h,\bar{a}\> }\iota
(a)$.

The untwisted space associated with $L$ is defined to be
\begin{eqnarray*}
V_{L}={\C}\{L\}\otimes _{{\C}}M(1)\simeq {\C}[L]\otimes S(\hat{\h}
^{-})\;\;\mbox{(linearly)}.
\end{eqnarray*} Then
$\hat{L},\hat{{\h}},z^{h}\;(h\in {\h})$ act naturally on $V_{L}$
by acting on either ${\C}\{L\}$ or $M(1)$ as indicated above.
Define ${\1}= \iota ( e_0) \in V_L$.
We use a normal ordering procedure, indicated by open colons,
which signify that in the enclosed expression, all creation
operators $h(n)$ $(n<0)$,$a\in \hat{L}$ are to be placed to the
left of all annihilation operators $h(n),z^{h}\;(h\in {\h},n\ge
0)$. For $a \in \hat{L}$, set
\begin{eqnarray*}
Y(\iota (a),z)= : e^{\int
(\bar{a}(z)-\bar{a}(0)z^{-1})}az^{\bar{a}}:.
\end{eqnarray*}
Let $a\in \hat{L};\;h_{1},\cdots,h_{k}\in
{\h};n_{1},\cdots,n_{k}\in {\Z}\;(n_{i}> 0)$. Set
\begin{eqnarray*}
v= \iota (a)\otimes h_{1}(-n_{1})\cdots h_{k}(-n_{k})\in
V_{L}.\end{eqnarray*}

 Define vertex operator $Y(v,z)$ with
\bea \label{defvertex}  :\left({1\over (n_{1}-1)!}({d\over
dz})^{n_{1}-1}h_{1}(z)\right)\cdots \left({1\over
(n_{k}-1)!}({d\over dz})^{n_{k}-1}h_{k}(z)\right)Y(\iota (a),z): .
\eea
 This gives us a well-defined linear map
\begin{eqnarray*}
Y(\cdot,z):& &V_{L}\rightarrow
(\mbox{End}V_{L})[[z,z^{-1}]]\nonumber\\ & &v\mapsto Y(v,z)=\sum
_{n\in {\Z}}v_{n}z^{-n-1},\;(v_{n}\in {\rm End}V_{L}).
\end{eqnarray*}

Let $\{\; h_{i}\;|\;i=1,\cdots,d\}$ be an orthonormal basis of
${\h}$ and set
\begin{eqnarray*}
\omega ={1\over 2}\sum _{i=1}^{d} h_{i}(-1) h_{i}(-1)\in V_{L}.
\end{eqnarray*}
Then $Y(\omega,z)=\sum_{n\in {\Z}}L_n z^{-n-2}$ gives rise to a
representation of the Virasoro algebra on $V_{L}$ with the central
charged $d$ and
\begin{eqnarray}  \label{vir.rel}
& &L_0\left(\iota(a)\otimes h_{1}(-n_{1})\cdots
h_{n}(-n_{k})\right)\nonumber \\ &=&\left({1\over 2}\<
\bar{a},\bar{a}\>+n_{1}+\cdots+n_{k}\right) \left(\iota(a)\otimes
h_{1}(-n_{1})\cdots h_{k}(-n_{k})\right).
\end{eqnarray}

The following theorem was proved in \cite{DL} and \cite{K}.
\begin{theorem}
$(V_L, {\1}, L_{-1}, Y)$ is vertex superalgebra.
\end{theorem}

Vertex algebra $M(1)$ can be treated as a subalgebra of $V_L$.

 Define the Schur polynomials $p_{r}(x_{1},x_{2},\cdots)$ $(r\in
{\Z}_{+})$ in variables $x_{1},x_{2},\cdots$ by the following
equation:
\begin{eqnarray}\label{eschurd}
\exp \left(\sum_{n= 1}^{\infty}\frac{x_{n}}{n}y^{n}\right)
=\sum_{r=0}^{\infty}p_{r}(x_1,x_2,\cdots)y^{r}.
\end{eqnarray}
For any monomial $x_{1}^{n_{1}}x_{2}^{n_{2}}\cdots x_{r}^{n_{r}}$
we have an element $h(-1)^{n_{1}}h(-2)^{n_{2}}\cdots
h(-r)^{n_{r}}{\1}$ in both $M(1)$ and $V_{L}$ for $h\in{\h}.$
 Then for any polynomial $f(x_{1},x_{2}, \cdots)$, \\ $f(h(-1),
h(-2),\cdots){\1}$ is a well-defined element in $M(1)$ and
$V_{L}$.  In particular, $p_{r}(h(-1),h(-2),\cdots){\1}$ for $r\in
{\Z}_{+}$ are elements of $M(1)$ and $V_{L}$.

Suppose $a,b\in \hat{L}$ such that $\bar{a}=\alpha,\bar{b}=\beta$.
Then
\begin{eqnarray}\label{erelation}
Y(\iota(a),z)\iota(b)&=&z^{\<\alpha,\beta\>}\exp\left(\sum_{n=1}^{\infty}
\frac{\alpha(-n)}{n}z^{n}\right)\iota(ab)\nonumber\\
&=&\sum_{r=0}^{\infty}p_{r}(\alpha(-1),\alpha(-2),\cdots)\iota(ab)
z^{r+\<\alpha,\beta\>}.
\end{eqnarray}
Thus
\begin{eqnarray}\label{eab1}
\iota(a)_{i}\iota(b)=0\;\;\;\mbox{ for }i\ge -\<\alpha,\beta\>.
\end{eqnarray}
Especially, if $\<\alpha,\beta\>\ge 0$, we have
$\iota(a)_{i}\iota(b)=0$ for all $i\in {\Z}_{+}$, and if
$\<\alpha,\beta\>=-n<0$, we get
\begin{eqnarray}\label{eab}
\iota(a)_{i-1}\iota(b)=p_{n-i}(\alpha(-1),\alpha(-2),\cdots)\iota(ab)
\;\;\;\mbox{ for }i\in {\Z}_{+}.
\end{eqnarray}

We will need one structural result on Heisenberg vertex algebras.

Element $L_0$ of the Virasoro algebra defines a
${\Z}_+$--graduation on vertex algebra $M(1)= \oplus_{n \in {\Z}_+
} M(1)_n$.
 Let $v_{\l}$ be the highest weight vector in the
$M(1)$--module $M(1, \l)$.

 The following lemma can be proved by using standard calculation in the
 Heisenberg vertex algebras.

\begin{lemma} \label{standard}
Let $h \in {\h}$, and $n_1, \dots, n_r \in {\N}$. Let $k= n_1 +
\cdots + n_r$. Let $u= h (-n_1) \cdots h (-n_r) {\1}$, and
$ Y( u, z) = \sum_{n \in {\Z} } u_{n} z^{-n-1}$.
Then $u \in M(1)_k$, and we have
\bea (1) && u_n v_{\l} = 0 \quad\mbox{for} \ n> k-1, \nonumber
\\
(2) && u_{k-1} v_{\l} = (-1) ^{n_1 + \cdots + n_r +r}( \la \l, h
\ra ) ^{r} v_{\l}. \nonumber \eea
\end{lemma}

\begin{proposition} \label{schur1}
Let $h \in {\h}$, and $r\in {\Z}_+$.  Let $u =
p_{r}(h(-1),h(-2),\cdots){\1} $. Set $ Y( u,z) = \sum _{n \in {\Z}
} u_n z^{-n-1}$. Then $u \in M_r$, and we have
\bea (1) && u_n v_{\l} = 0 \quad \mbox{for} \ n>r-1, \nonumber \\
(2)&& u_{r-1} v_{\l} = { \la \l, h \ra \choose r} v_{\l}. \eea
\end{proposition}
{\em Proof.} Since $u \in M(1)_r$, we have that $u_n v_{\l} = 0$
for $n > r$. Set $x=\la \l , h \ra$. Using Lemma \ref{standard}
one can easily see that
$$ u_{r-1} v_{l}= p_r (x, -x, x,-x,\cdots )  v_{\l}.$$ Since
$$ \exp \left( \sum_{n=1} ^{\infty} \frac{(-1) ^{n-1} x}{n} y^{n}
\right) = ( 1+y) ^{x} = \sum_{r \ge 0} {x \choose r} y^{r},$$
we have that $ p_r (x, -x, x,-x,\cdots ) = {x \choose r}$, and we
conclude that $u_{r-1} v_{\l}= {x \choose r}$. \qed

\begin{remark} Proposition \ref{schur1} can be also proved by
using Zhu's algebra theory (see \cite{DLM}, Section 3.)
\end{remark}

\section{Representations of the vertex algebra $W_{1+\infty, -N}$ }
\label{wsec}

 In this section we will construct $2N$--dimensional
family of irreducible modules for the vertex algebra $W_{1+\infty,
-N}$. We will prove that $W_{1+\infty, -N}$ can be realized as a
vertex subalgebra of Heisenberg vertex algebra $M(1)$.

First we will consider the following lattice:
\bea
&& L=\bigoplus_{i= 1} ^{N} \left( {\Z}  {\a}_i + {\Z} {\b} _i
\right), \nonumber \quad \mbox{where}\\
&& \la {\a}_i , {\a}_j \ra = \delta_{i,j}, \quad \la {\b}_i ,
{\b}_j \ra = - \delta_{i,j}, \quad \la {\a}_i, {\b}_j \ra = \la
{\b}_j , {\a}_i \ra = 0, \nonumber \eea
for every $i,j \in \{1,\cdots,N\}$.

 Let ${\h}$, $M(1)$ and $V_L$
be defined as in the Section \ref{stwist}. Then ${\h}$ is a
Heisenberg algebra, $\mbox{dim} {\h} = 2 N$, $M(1)$ is a
Heisenberg vertex algebra, and $V_L$ is a lattice vertex
superalgebra.

 For every $ i\in \{1, \cdots,
N\}$, let $a^{i}, b^{i} \in {\hat L}$ such that
$${\overline{a^{i} } } = \a_i + \b_i, \quad {\overline{b^{i} }}=-
(\a_i + \b_i).$$
 Then we define $e^{i},
f^{i} \in V_L$ with
$ e^{i} = \iota( a^{i} ), f^{i} = \iota( b^{i})$.

Set $Y( e^{i}, z) = \sum_{n \in {\Z} } e^{i} _n z^{-n-1}$, and
$Y( f^{i}, z) = \sum_{n \in {\Z} } f^{i} _n z^{-n-1}$.

Recall the definition of Schur polynomial  $p_l ( h(-1),h(-2),
 \dots)$ ($h\in {\h}$ ) from Section \ref{stwist}.
Set $p_l (h) := p_l ( h(-1),h(-2), \dots) \in S({\hat{\h}} ^{}-)$.
Let $\gamma_i = \a_i + \b_i$. Then $ \la \gamma_i, \gamma_i \ra
=0$.

 Now, relations (\ref{eab1}) and (\ref{eab}) imply the following
 lemma.

\begin{lemma} \label{fms1}
For all $i,j \in \{1, \cdots, N\}$, $k,n \in {\Z}$
we have :
\bea (1) && e^{i} _n  e^{j} = f^{i} _n f^{j}= 0, \quad \mbox{for}
\ \ n \ge 0, \nonumber \\
(2) && e^{i} _{-1}  f^{j} = \delta_{i,j} {\1}, \quad  f^{i} _{-1}
e^{j}= \delta_{i,j} {\1}, \nonumber \\
(3) && e^{i} _{-l-1}  f^{i} = p_l( \gi) {\1}, \quad  f^{i} _{-l-1}
e^{i}= p_l (-\gi) {\1}, \nonumber \\
(4)&& [ \a_i (k), e^{j} _n ]= \delta_{i,j} e^{j} _{n+k}, \quad [
\a_i (k), f^{j} _n ]= -\delta_{i,j} f^{j} _{n+k}.  \eea
\end{lemma}

  Define
\bea &&A^{i} ( z) = \sum_{n \in {\Z} } A^{i} (n)
z^{-n-1}=Y(e^{i}\otimes \alpha_i (-1), z ) =  \ :\alpha_i (z)
Y(e^{i} , z) :, \nonumber
\\
&& {\bA}^{i} (z) = \sum_{n \in \Z} {\bA} ^{i} (n) z^{-n}=
Y(f^{i} ,z) . \nonumber  \eea

\begin{lemma}  \label{fms2} We have
  \bea
  && [A^{i} (n), A ^{j} (m)]= [{\bA} ^{i} (n), {\bA} ^{j} (m)]= 0,
  \nonumber \\
  && [{\bA} ^{i} (n), A^{j} (m) ] = \delta_{i,j} \delta_{m+m,0},
  \nonumber \eea
  i.e. the components of the fields
  $A^{i}(z), {\bA} ^{i} (z)$, $1\le i \le N$,  span an associative algebra isomorphic to the
  Weyl algebra $W_N$ with $C=1$.
 \end{lemma}
 {\em Proof.}
Using (\ref{5.2}) we have that
$$ A^{i}(m) = \sum_{ n < 0} \a_i (n) e^{i }_{ m-n-1} + \sum_{n \ge
0} e^{i}_{m-n-1} \a_i (n) $$ for $1\le i \le N$.
Using Lemma \ref{fms1} we get
\bea
&& A^{i} (0) f^{j}= -\delta_{i,j} f^{j}, \quad A^{i} (m) f^{j} = 0
\quad \mbox{for} \ m > 0, \nonumber \\
&& A^{i} (m) A^{j} (-1){\1} = 0, \quad f^{i} _m f^{j} = 0, \quad
\mbox{for} \ m\ge 0. \nonumber \eea
Now, the statment of the lemma follows from the commutator
formulae  (\ref{comut}). \qed

  From Lemma \ref{fms2} follows the following result.

\begin{proposition} \label{wizom}The subalgebra of the vertex superalgebra $V_L$
generated by the fields $A^i (z)$, ${\bA} ^{i} (z) $ $i=1, \dots,
N$, is isomorphic to the vertex algebra $V_N$. (Under this
isomorphism the fields $a_i (z) $ corresponds to $A^{i} (z)$, and
$a^{*} _i (z)$ to $\bar{A} ^{i}(z)$.)
\end{proposition}
{\em Proof.} From the Lemma \ref{fms2} follows that $V_L$ is a
module for the Weyl algebra ${\Wi}$. The subalgebra of $V_L$
generated by the fields $A^i (z)$ and ${\bA} ^{i} (z) $ is exactly
the ${\Wi}$ submodule ${\Wi} . {\1}$ generated by ${\1}$. Then we
have
$$ A^{i} (m) {\1} = {\bA} ^{i} (n) {\1} = 0, \quad \mbox{for} \
m\ge0, \ n > 0.$$ Then the  fact that $V_N$ is an irreducible
${\Wi}$--module implies
 that $V_N \cong {\Wi}. {\1} $. \qed

\begin{remark} In the physical literature the vertex algebra $V_N$
is known as $\beta \gamma$ system, and the lattice construction of
$\beta \gamma$ system is known as Friedan-Martinec-Shenker
bosonization (cf. \cite{FMS}, \cite{W1}). \end{remark}

 With the respect to Proposition
\ref{wizom} we can identify vertex algebra $V_N$ with the
subalgebra of $V_L$ generated by the fields $A^{i}(z)$, $\bar{A}
^{i}(z)$. Since $W_{1 +\infty,-N}$ is a vertex subalgebra of $V_N$
generated by the fields
$$ J^{k} (z) = - \sum_{ i=1} ^{N} : A ^{i} (z)\partial_z ^{k}
{\bA} ^{i} (z):$$
we have that $W_{1 +\infty,-N}$ is also a vertex subalgebra of
$V_L$. For every $i \in \{1, \cdots, N\}$, let $U^{i}_{k}= A^{i}
(-1)f^{i} _{-k-1} {\1} \in V_L$, and $U_{k}= \sum_{i=1} ^{N} U^{i}
_{k}$ . We have
$$ J^{k} (z) = -k! \sum_{i=1} ^{N } Y ( A^{i} (-1) {\bA}^{i} (-k)
{\1},z ) =-k! Y(U_k,z).$$

\begin{lemma} \label{ul1} We have
\bea
&&U^{i} _{k} = A^{i}(-1) {\bA}^{i} (-k) {\1}= \a_i (-1) p_k (-\gi)
{\1} + p_{k+1} (-\gi)  {\1}. \nonumber \eea In particular,
$U^{i}_{k} \in M(1)$.
\end{lemma}
{\em Proof.} Since $[A^{i} (-1), {\bA}^{i} (-k)] = 0$, we have
that
$$
 U^{i} _k = {\bA}^{i} (-k)A^{i}(-1) {\1} = f^{i} _{-k-1} \a_i
(-1) e^{i}= \a_i (-1) f^{i}_{-k-1} e^{i} + f^{i} _{-k-2} e^{i}.$$
Then Lemma \ref{fms1} implies that
$$ U^{i} _{k} = \a_i (-1) p_k (-\gi) {\1} + p_{k+1} (-\gi) {\1},
$$ and lemma holds. \qed

\begin{theorem} Vertex algebra $W_{1+ \infty, -N}$ is a subalgebra of the  vertex algebra $M(1)$.
\end{theorem}
{\em Proof.} From Lemma \ref{ul1} follows that $U_{k}= \sum_{i=1}
^{N} U^{i} _{k} \in M(1)$, for every $k \in {\Z}_+$. Since $J^{k}
(z) = - k! Y(U_{k},z)$, $k \in {\Z}_+$, generate $W_{1+ \infty,
-N}$, we have that $W_{1+ \infty, -N}$ is a subalgebra of the
vertex algebra $M(1)$. \qed

\begin{lemma} \label{hw}
For every $\l \in {\h}$, let $v_{\l}$ be the highest weight vector
in $M(1,\l)$. Then we have
$ J ^{k}(n)  v_{\l} = 0$ for $ n > 0$, and
$$ J ^{k} (0) v_ {\l} = - k! \sum_{i=1} ^{N} \left( {-\la \l , \gi
\ra \choose {k+1} }  + \la \l, \a_i \ra {-\la \l, \gi \ra \choose
k} \right) v_{\l} . $$
\end{lemma}
{\em Proof.} Set $Y( U^{i} _k, z) = \sum_{n \in {\Z} } U^{i} _k
(n) z^{-n-k-1}$. Proposition \ref{schur1}, and Lemma \ref{ul1}
implies that
\bea \label{ul2}&&  U^{i} _k (n) v_{\l} = 0 \quad \mbox{for} \ n
>0, \nonumber
\\ && U^{i} _k (0) v_{\l} = \la \l, \a_i \ra  {-\la \l , \gi \ra \choose k} +
{-\la \l, \gi \ra  \choose k+1}. \nonumber \eea
Since $ J^{k} (n) = - k!\sum_{i=1} ^{N} U^{i} _k (n)$, we have
that the
 lemma holds. \qed

  Since $W_{1+ \infty, -N}$ is a subalgebra of VOA $M(1)$,
then for every $\l \in {\h}$, $M(1, \l)$ is a $W_{1+ \infty,
-N}$--module. Let $V(\l,-N)$ be the irreducible $W_{1+ \infty,
-N}$--subquotient of $M(1, \l)$ generated by the vector $v_{\l}$.

\begin{theorem} \label{main}
For every $\l \in {\h} $ $ {V} ({\l},-N)$ is the irreducible
  module for the vertex algebra $W_{1 +\infty, - N}$.
As a $\hD$--module, $V(\l,-N)$ is a irreducible quasi-finite
highest weight module, and the corresponding generating series is
\bea \label{gener}
\Delta_{\lambda} (x) =-\sum_{i=1} ^{N}
\left(\frac{e^{ s_i x} -1 }{e^{x} -1} + t_i e^{s_i x} \right) ,
\eea
where
\bea \label{defl}s_i = -\la \l, \a_i + \b_i \ra, \quad t_i = \la
\l, \a_i \ra, \quad i=1,\dots, N. \eea
\end{theorem}
{\em Proof.} Lemma \ref{hw} implies that ${V}({\l},-N)$ is a
highest weight $\hD$--module. It remains to prove that ${
V}({\l},-N)$ is a quasi-finite $\hD$--module. In order to prove
this we have to identify the generating series $\Delta_{\l} (x)$.

Using the relation (\ref{veza}), it is straightforward to see the
following formulae.
\bea  \Delta_{\l} (x) = \sum_{k=0} ^{\infty} \frac{x^{k}}{k!}
L^{k} (0) v_{\l} = \sum_{k=0} ^{\infty} \frac{(e^{x} -1)^{k}}{k!}
J^{k} (0) v_{\l}. \nonumber \eea

Let $s_i, t_i$ be defined  with (\ref{defl}). Then Lemma \ref{hw}
implies
\bea \Delta_{\l} (x) = &&-\sum_{i=1} ^{N}\sum_{k=0} ^{\infty}
\left( {s_i \choose {k+1} }  + t_i {s_i \choose k} \right) (e^{x}
-1)^{k} \nonumber
\\
= && -\sum_{i=1} ^{N} \left(\frac{e^{s_i x} -1 }{e^{x} -1} + t_i
e^{s_i x} \right) \nonumber
\\
=&& \frac{\Phi (x) }{e^{x} -1}, \nonumber \eea
where
$$ \Phi(x) -N = \sum_{i=1} ^{N} \left( (t_i -1 ) e^{s_i x} -t_i
e^{(s_i +1) x}\right).$$ This implies that ${V}(\l,-N)$ is a
quasi-finite $\hD$--module. \qed

\begin{remark}
Theorem \ref{main} gives the existence of $2 N$--dimensional
family of irreducible $W_{1+\infty,-N}$--modules. If we take in
(\ref{gener})  $t_i=0$ for every $i =1, \dots, N$, we get exactly
$W_{1+\infty,-N}$--modules constructed in \cite{KR2}.
\end{remark}

 We
have the following conjecture.

\begin{conjecture} \label{slutnja}
The set $V(\l,-N)$, $\l \in {\h}$, lists all the irreducible
modules for the vertex algebra $W_{1+\infty,-N}$.
\end{conjecture}

\begin{remark} In Section \ref{sec-1} we will see that the Conjecture \ref{slutnja}
is true for $c=-1$. \end{remark}

\section{ The case of $c=-1$ } \label{sec-1}

In this section we will compare our results with the results from
\cite{W1}, \cite{W2}. In \cite{W1}, Wang proved that the vertex
algebra $W_{1+\infty,-1}$ is isomorphic to the tensor product
$W_{3,-2} \otimes H_0$, where $W_{3,-2}$ is a simple vertex
algebra associated to $W_3$--algebra with $c=-2$, and $H_0$ is a
Heisenberg vertex algebra. Moreover, in \cite{W2} Wang classified
all the irreducible modules for $W_{3,-2}$ and $W_{1+\infty,-1}$.
The methods used in \cite{W1}, \cite{W2} didn't imply the
identification of $W_{1+\infty,-1}$--modules as a highest weight
$\hD$--modules (see Section 5 in \cite{W2}). Our approach gives an
explicit identification of two-dimensional family
$W_{1+\infty,-1}$--modules in terms of highest weights.

 Let $N=1$, and set $\a= \a_1$, $\b = \b_1$.
For $\l \in {\h}$ set $\l_{\a} = \la \l, \a \ra$, $\l_{\b} = \la
\l , \b \ra$. Let $M_{\a}(1, \l_{\b})$  (resp. $M_{\b} (1,
\l_{\b})$ ) be the submodules of $M(1,\l)$ generated by the
highest weight vector $v_{\l}$ and $\a (n)$ (resp. $\b (n)$).  Set
$M_{\a} (1) = M_{\a}(1,0)$, $M_{\b} (1) = M_{\b}(1,0)$. Then
\bea \label{m1ab} M(1) = M_{\a} (1)\otimes M_{\b} (1), \quad
M(1,\l) = M_{\a} (1, \l_{\a} )\otimes M_{\b} (1, \l _{\b}). \eea

As in \cite{W2} we define

\bea
&& T(z) = \frac{1}{2} \left( :\a (z) ^{2} : + \partial\a(z)
\right), \nonumber \\
&& W(z) = \frac{1}{12} \left( 4 : \a (z) ^{3} : + 6 : \a(z)
\partial\a(z) : + \partial^{2} \a(z) \right). \nonumber
\eea

\begin{theorem} \label{wang} \cite{W1}, \cite{W2}

 \item[(1)] The fields $T(z)$, $W(z)$ span a subalgebra of $M_{\a}(1)$ isomorphic
to $W_{3,-2}$, and  $W_{1+\infty,-1}\cong W_{3,-2} \otimes M_{\b}
(1)$.
\item[(2)] Let ${\cal V}_r$ be the irreducible subquotient of
$W_{3,-2}$--module $M_{\a}(1,r)$. Then ${\cal V}_r$, $r \in {\C}$,
gives all the irreducible $W_{3,-2}$--modules.
\end{theorem}

Recall the definition of $W_{1+\infty,-1}$--modules $V(\l,-1)$
from Section \ref{wsec}. Then we have the following consequence of
Theorem \ref{main} and Theorem \ref{wang}.

\begin{corollary} We have
\item[(1)] $ {V}(\l,-1) \ = {\cal V}_{\l_{\a} } \otimes M_{\b} (1,
\l_{\b} )$ \ \ for every $\l \in {\h}$.
\item[(2)]The set $V(\l,-1)$, $\l \in {\h}$, gives all irreducible
$W_{1+\infty,-1}$--modules.
\item[(3)] $$\Delta_{\l} (x) =-\frac{e ^{-(\l_ {\a} + \l_{\b}) x
}-1}{e^{x}-1}- \l_{\a} e^{- (\l_ {\a} + \l_{\b}) x}.$$
\end{corollary}
{\em Proof.} (1) follows from the definition of ${\cal V}_{\a}$
and (\ref{m1ab}).  Then theorem \ref{wang} implies that ${\cal
V}_r \otimes M_{\b} (1,s)$, $r,s \in {\C}$, are all irreducible
$W_{1+\infty,-1}$--modules. Since
 $ {V}(\l,-1) \ = {\cal V}_{\l_{\a} } \otimes M_{\b} (1,
\l_{\b} )$, we see that $V(\l,-1)$, $\l \in \h$ gives all
irreducible $W_{1+\infty,-1}$--modules, and we get (2). The
statement (3) follows from the Theorem \ref{main}. \qed

Department of Mathematics, University of Zagreb, Bijeni\v{c}ka 30,
10000 Zagreb, Croatia

E-mail address: adamovic@cromath.math.hr
\end{document}